\documentclass{elsart3-1}

\usepackage{amssymb}

\usepackage[english,francais]{babel}

\newtheorem{theorem}{Theorem}[section]
\newtheorem{lemma}[theorem]{Lemma}
\newtheorem{proposition}[theorem]{Proposition}
\newtheorem{corollary}[theorem]{Corollary}

\renewcommand{\theequation}{\arabic{equation}}
\def\og{\leavevmode\raise.3ex\hbox{$\scriptscriptstyle\langle\!\langle$~}}
\def\fg{\leavevmode\raise.3ex\hbox{~$\!\scriptscriptstyle\,\rangle\!\rangle$}}
\def\EMdash{\leavevmode\hbox to 7.5mm{\vrule height .63ex depth -.59ex width 5.4mm\hfill}}

\let\epsilon=\varepsilon

\begin{document}

\begin{frontmatter}




%
\selectlanguage{francais}
\title{\textbf{Anneaux d'endomorphismes et classe d'isog\'enies de modules de Drinfeld de
 rang $2$ sur un corps fini}}

\vspace{-2.6cm} \selectlanguage{english}
\title{\textbf{Endomorphism Rings and Isogenies Classes for Drinfeld $A$-Modules of
Rank $2$ over Finite Fields}}



\author[a]{Mohamed-Saadbouh MOHAMED-AHMED}
\ead{mohamed-saadbouh.mohamed-ahmed@univ-lemans.fr}

\address[a]{D\'epartement de Math\'ematiques, Universit\'e du Maine, Avenue Olivier Messiaen,
72085 Le Mans Cedex 9, France}

\selectlanguage{francais}
\begin{abstract}
Soit $\Phi $ un $\mathbf{F}_{q}[T]$-module de Drinfeld de rang
$2$, sur un corps fini $L$, une extension de degr\'e $n$ d'un corps
fini $\mathbf{F}_{q}$. Nous sp\'ecifions les conditons de
maximalit\'e et de non maximalit\'e pour l'anneau d'endomorphismes
End$_{L}\Phi$ en tant que $\mathbf{F}_{q}[T]$-ordre dans
l'anneau de division End$_{L}\Phi \otimes _{\mathbf{F}_{q}[T]}\mathbf{F}%
_{q}(T)$, on s'intr\'essera ensuite aux polynôme caract\'eristique et par son
interm\'ediaire on calculera le nombre de classes d'iog\'enies.

{\it Pour citer cet article:
Mohamed-saadbouh.Mohamed-Ahmed, C. R. Acad. Sci. Paris,
 Ser. I ... (...).}

\vskip 0.5\baselineskip \selectlanguage{english}
\noindent{\bf Abstract} \vskip 0.5\baselineskip \noindent We
discuss many analogy points with the elliptic curves. more
precisely, we study the characteristic polynomial of a Drinfeld
module of rank $2$ and use it to calculate the number of isogeny
classes for such modules.

{\it To cite this article: Mohamed-Saadbouh Mohamed-Ahmed
, C. R. Acad. Sci. Paris, Ser. I ... (...).}
\end{abstract}

\end{frontmatter}

\section{Introduction}

 let $K$ a no empty global field of characteristic $p$ (
 namely a rational functions field of one indeterminate over a
 finite field ) together with a constant field, the finite field
 $\mathbf{F}_{q}$ with $p^{s}$ elements. We fix one place of $K$,
 denoted by $\infty ,$ and call $A$  the ring of regular elements
 away from the place $\infty $. Let $L$ be a commutator field of
 characteristic $p$,  $\gamma :A\rightarrow L$ be a ring
 $A$-homomorphism. The kernel of this $A$-homomorphism is denoted by $P.$
 We put $m$ =$[L,$ $A/P]$, the extension degrees of $L$ over $A/P$.

 We denote by $L\{\tau \}$ the Ore polynomial ring, namely, the
 polynomial ring of $\tau $, where $\tau $ is the Frobenius of
 $\mathbf{F}_{q}$ with the usual addition and where
 the product is given by the commutation rule: for every $%
 \lambda \in L\mathit{,}$ we have $\tau \lambda =\lambda ^{q}\tau
 $. A Drinfeld $A$-module  $\Phi :A \rightarrow
 L\{\tau \}$ is a non trivial ring homomorphism and a non trivial embedding
 of $A $ into $L\{\tau \}$ different from $\gamma $.
 This homomorphism $\Phi $, once defined, define an $A$-module structure over the $A$-field $L$, noted \ $%
 L^{\Phi }$, where the name of a Drinfeld $A$-module for a
 homomorphism $\Phi $. This structure of $A$-module depends on
 $\Phi $ and, especially, on his rank.

   Let $\chi $ be the Euler-Poincaré characteristic ( i.e. it is an
   ideal from $A$). So we can speak about the ideal $\chi (L^{\Phi
   })$, denoted henceforth by $\chi _{\Phi }$, which is by definition
   a divisor of $A$, corresponding for the elliptic curves to a
   number of points of the variety over their basic
   field. In this paper, we will work on the special case $K=$ $\mathbf{F}%
   _{q}(T)$, $\ A=\mathbf{F}_{q}[T]$. Let $P_{\Phi }(X)$ be the
   characteristic polynomial of the $A$-module $\Phi $, which is also
   a characteristic polynomial of the Frobenius $F$ of $L$.\ We can
   prove that this polynomial can be given as:
   $P_{\Phi }(X)=$ $X^{2}-cX+\mu P^{m},$ such that $\mu \in \mathbf{F}%
   _{q}^{\ast } $, and $c$ $\in A,$ where $\ \deg c\leq
   \frac{m.d}{2}$ by the Hasse-Weil analogue in this case. We will be
   interested in the endomorphism ring and the isogeny classes numers  of a
   Drinfeld $A$-module of  rank 2. for more information see
   \cite{Angeles}, \cite{B.Angles},  \cite{Gosse}, and
   \cite{Drinfeld2}.

\subsection{The endomorphism ring}
 The Drinfeld $A$-module
 of rank $2$ is of the form $\ \Phi (T)=a_{1}+a_{2}\tau +a_{3}\tau
 ^{2}$, where $\ a_{i}\in L$, $1\leq i\leq 2$ and $a_{3}\in L^{\ast
 }$. Let $\Phi $ and $\Psi $ be two Drinfeld modules over an $A$-field
     \ $L$. A morphism from $\Phi $ to $\Psi $ over $L$ is an element $ p(\tau )\in L%
     \mathit{\{\tau \}}$ such that
     $p\Phi _{a}=\Psi _{a} p\quad \mbox{ for all }
     a\in A$.
   A non-zero morphism is called an isogeny. We note that this is
   possible only between two Drinfeld modules with the same rank.
   The set of all morphisms forms an $A-$module denoted by Hom$_{E}(\Phi ,\Psi ).$

    In particular, if $\Phi $ =$\Psi $ the $L$-endomorphism ring (End$%
    _{L}\Phi =$Hom$_{L}(\Phi ,\Phi $) is a subring of \ $L\{\tau \}$ and an $A$%
    -module contained in $\Phi (A)$. Let $F$ be the Frobenius of $L$ we have : $\Phi $ $(A)\subset $ \
    End$_{L}\Phi$ and $\ F\in $ End$_{L}\Phi $.

Let $\overline{L}$ be a fix algebraic closure of $L$,

$\Phi _{a}(\overline{L}):=\Phi \lbrack a](\overline{L})=\{x\in \overline{L}%
,\Phi _{a}(x)=0\}$, and  $\Phi _{P}(\overline{L})=\cap _{a\in P}\Phi _{a}(%
                                                  \overline{L}).$
     We say that $\Phi $ is supersingular if and only if the $A$-module constituted by a $P$%
     -division points $\Phi _{P}(\overline{L})$ is trivial,
     otherwise  $\Phi $ is said a ordinary module.

\begin{proposition}
Let $P_{\Phi}(X)=X^{2}-cX + \mu P^{m}$ be the characteristic
polynomial of the Frobenius $F$ of a finite field $L$ and let
$\Delta = c^{2}-4\mu P^{m}$ be the discriminat of $P_{\Phi}$,
and $O_{K(F)}$ the maximal $A$-order of the algebra $K(F)$.
\begin{enumerate}
\item For every $g \in A $ such that $\Delta = g^{2}.\omega $,
there exists a Drinfeld $A$-module $\Phi $ over $L$ of rank $2$ such
that $O_{K(F)}=A[ \sqrt{ \omega }]$ and:
$End_{L} \Phi =A+ g.A[\sqrt{\omega }]. $
\item If there is no polynomial $g$ of $A$ such that $g^{2}$
divide $ \Delta ,$ then there exists an ordinary Drinfeld
$A$-module $\Phi $ over $L$ of rank $2$ such that
$End_{L}\Phi =O_{K(F)}.$
\end{enumerate}
\end{proposition}

\subsection{ Isogeny classes}

Let $\overline{K}$ be an algebraic closure of $K$ and let $\infty
$ be a place of $K$ which divides $\frac{1}{T}.$ Let us put
$K_{\infty }=F_{q}((\frac{ 1}{T}))$ and denote by
$\mathbb{C}_{\infty }$ the completude of the algebraic closure of
$K_{\infty }$. We fix an embedding $\overline{K}$
$\hookrightarrow \mathbb{C}_{\infty }$. For every $\alpha \in
\mathbb{C}_{\infty }$, we denote by $\mid \alpha \mid _{\infty }$
the normalized valuation of $\alpha $ ( $\mid \frac{1}{T}\mid _{\infty }=%
\frac{1}{q}$).

Let $\theta \in $ $\overline{K}$, we say that $\theta $ is an
ordinary number if:

\begin{enumerate}
\item $\theta $ is integral over $A$;

\item $\mid \theta \mid _{\infty }=q^{\frac{md}{2}};$

\item $K$($\theta )/K$ is imaginary and $[K$($\theta
),K\mathbf{]=}2;$

\item there is only one place of $K$($\theta )$ which divides $\theta $ and
Tr$ _{K(\theta )/K}(\theta )\neq 0(P).$
\end{enumerate}
We say that $\theta $ is an ordinary Weil number if
$\theta^{\sigma }$ is an ordinary number for all $\sigma \in $
col($\overline{K}\mathbf{/}K$). We denote by W$^{ord}$ the set of
conjugancy class of ordinary Weil numbers of rank 2. We have the
important result, for proof  see \cite{YU}:

\begin{theorem}
There exists a bijection between W$^{ord}$ and isogeny classes of
ordinary Drinfeld $A$-modules of rank $2$ defined over $L$.
\end{theorem}
Let $\theta $ be an ordinary Weil number. We put
$P(x)=Hr(\theta, K\mathbf{;}x\mathbf{)}.$
By using (1), (2), (3) and (4) we have
$P(x)=x^{2}-cx+\mu P^{m}$,
where $\mu \in \mathbf{F}_{q}^{\ast }$ and $c$ $\in A$ but $c\neq 0(P)$ and also deg$%
_{T}$ $c\leq \frac{md}{2}.$ Let us put
$$
\Gamma =\left \{c \in A \quad \mbox{ such that } c\neq 0(P) \mbox{
and } {\mathrm deg}_{T} c\leq \frac{md}{2}\right \}.$$
We need the following lemma.

\begin{lemma}
For $\mu \in \mathbf{F}_{q}^{\ast }$, $c \in \Gamma $ denote by
$E$ the filed of decomposition of $P(x)=x^{2}-cx+\mu P^{m}$ over
$K$. Let $\theta $ be a root of $P(x)$. Then $\theta $ verifies
(1), (2), (3) and (4) together with $[K(\theta ),K\mathbf{ ]=}2$.
\end{lemma}

\begin{corollary}
1) Let $\mu \in \mathbf{F}_{q}^{\ast }$ and $c$ $\in \Gamma $ and
let $ \theta $ be a root of $x^{2}-cx+\mu P^{m}.$ Then $\theta $
is an ordinary Weil number if and only if $K(\theta )/K$ is
imaginary.

2) If $ md \equiv 1(2)$, then the roots of $x^{2}-cx+\mu P^{m}$ are
Weil numbers for all $\mu \in \mathbf{F}_{q}^{\ast }$ and for all
$c\in \Gamma.$
\end{corollary}
To simplify, let us suppose $p\neq 2$ and put $md\equiv 0(2).$
\begin{lemma}
Let $\mu \in \mathbf{F}_{q}^{\ast }$ and $c$ $\in \Gamma $ with
deg$_{T}$ $c\leq \frac{md}{2}$. Let $\theta $ be a root of
$x^{2}-cx+\mu P^{m}$.
Then $\theta $ is a Weil number if and only if $-\mu \notin (\mathbf{F}%
_{q}^{\ast })^{2}$.
\end{lemma}

\begin{lemma}
Let $\mu \in \mathbf{F}_{q}^{\ast }$ and $c$ $\in \Gamma $ with deg$_{T}$ $%
c= \frac{md}{2}.$ Denote by $c_{0}$ the term of higher degree of
$c$. We suppose that $c_{0}^{2}\neq -4\mu $. Let $\theta $ be a
root of $x^{2}-cx+\mu P^{m}$. Then $\theta $ is a Weil number if
and only if $x^{2}-c_{0}x+\mu $ is irreducible in
$\mathbf{F}_{q}[X]$.
\end{lemma}
 The roots of the characteristic polynomial are a Wiel Numers, so
 we need this result, for proof see  \cite{YU}:
\begin{proposition}
Let $\Phi $ be a Drinfeld $A$-module of rank $2$ over the finite
field $L= \mathbf{F}_{q^{n}}$ and let $P$ be the characteristic of
$L$. We put $m=[L:A/P]$ and $d=$deg $P$. The characteristic
polynomial $P_{\Phi }$ can take only the following forms:
\begin{enumerate}
\item In the case of ordinary Drinfeld $ A$-modules, we have
$P_{\Phi }(X)=X^{2}-cX+\mu P^{m}$, where $c^{2}-4\mu P^{m} $ is
imaginary, $c\in A$, $(c,P)=1$ and $\mu \in \mathbf{F}_{q}^{\ast
}.$

\item In the case for supersinglar $A$-modules, we distinguish
three cases:
\begin{enumerate}
\item If $m$ is odd, then $P_{\Phi }(X)=X^{2}+\mu P^{m}$, with
$\mu \in \mathbf{F}_{q}^{\ast }.$

\item If $m$ is even and $d=\deg P$ is odd, then $\ \ P_{\Phi
}(X)=X^{2}+c_{0}X+\mu P^{m},$ with $\mu \in \mathbf{F}_{q}^{\ast
}$ and $c_{0}\in \mathbf{F}_{q}$.

\item If $m$ is even, then $P_{\Phi }(X)=(X+\mu
P^{\frac{m}{2}})^{2}.$
\end{enumerate}
\end{enumerate}
\end{proposition}

 We can recapitulate all the cases above as follows :

 \begin{enumerate}
 \item For the ordinary case, the characteristic polynomial is of the form :%
 $P_{\Phi }(X)=X^{2}-cX+\mu P^{m} $,
 such that $2$ deg $c<\deg P.m$ or $2\deg c=$deg $P.m$ and $%
 X^{2}-a_{0}X+\mu $ is irreducible over $\mathbf{F}_{q^{n}}$ where
 $a_{0}$ is the coefficient of the greatest degree of $c$. For the
 supersingular case, we have the two following cases :

 \item The $\deg P$ is even or $-\mu \notin (\mathbf{F}_{q}^{\ast
 })^{2}$.

 \item the polynomial $X^{2}+c_{0}X+\mu $ is irreducible over
 $\mathbf{F}_{q}$.
 \end{enumerate}

We are in position to compute the number of characteristic
polynomials which corresponding to the number of isogeny classes, for proof see \cite{Gekeler}.
\begin{lemma}
$\#$\{Isogeny classes\}=$\#\{P_{\Phi }\}.$
\end{lemma}

So we can compute the cardinal of the isogeny classes of a
Drinfeld module of rank $2 $ as follows.
\begin{proposition}
Let $\Phi $ a Drinfeld $A$-module of rank $2$ over a finite field $%
L=F_{q^{n}}$ and  let $P$ be the $A$-characteristic of $L$. We
put $m=[L:A/P] $ and $d=$deg $P$ :
\begin{enumerate}
\item If $m$ and $d$ are both odd, then
$\#\{P_{\Phi }\}=(q-1)(q^{[\frac{m}{2}d]+1}-q^{[\frac{m-2}{2}d]+1}+1). $

\item If $m$ is even and $d$ is odd, then
$\#\{P_{\Phi
}\}=(q-1)[\frac{q-1}{2}q^{\frac{m}{2}d}-q^{\frac{m-2}{2}d+1}+q].  $

\item If $m$ and $d$ are both even, then
$\#\{P_{\Phi
}\}=(q-1)[\frac{q-1}{2}q^{\frac{m.}{2}d}-q^{\frac{m-2}{2}d}+1].  $
\end{enumerate}
\end{proposition}

\subsubsection{Euler-Poincare characteristic}

\bigskip Let $\Phi $ be a Drinfeld $A$-module of rank $2$ over a finite
field $L=\mathbf{F}_{q^{n}}$ and denote by $P_{\Phi }$ the
characteristic polynomial. Let $\chi _{\Phi }=(P_{\Phi }(1))$ This
Euler-Poincar\'{e} characteristic

We can have an expression for the cardinal of the set of
Euler-Poincare characteristic as follows.
\begin{proposition}
Let $\Phi $ be a Drinfeld A-module of rank $2$ over the finite field  $L=%
\mathbf{F}_{q^{n}}$, and let $P$ be the characteristic of $L$. We put $%
m=[L:A/P]$ and $d=$deg $P$. There exists $H,B\in L,$ such that
$\#\{\chi _{\Phi }\}=H+B,$
where $H$ and $B$ satisfies
$\#\{P_{\Phi }\}=(q-1)H+(q-2)B. $
\end{proposition}

The value of $\#\{\chi _{\Phi }\}$ can be deduced accordingly.

\begin{proposition}
Let $\Phi $ be a Drinfeld A-module of rank $2$ over a finite field $L=%
\mathbf{F}_{q^{n}}$ and let $P$ be the $A$-characteristic of $L$. We put $%
m=[L:A/P]$ and $d=$deg $P$. We have:

\begin{enumerate}
\item If $m$ and $d$ are both odd, then

$\#\{\chi _{\Phi }\}=\frac{q}{q-1}q^{[\frac{m}{2}d]+1}-\frac{q}{q-1}q^{[\frac{%
m-2}{2}d]+1}+1. $

\item If $m$ is even and $d$ is odd, then
$\#\{\chi _{\Phi }\}=\frac{q^{2}+1}{2q-2}q^{\frac{m}{2}d}-\frac{q}{q-1}q^{%
\frac{m-2}{2}d+1}+q. $

\item If $m$ and  $d$ are both even, then
$\#\{\chi _{\Phi }\}=\frac{q^{2}+1}{2q-2}q^{\frac{m}{2}d}-\frac{q}{q-1}q^{%
\frac{m-2}{2}d+1}+1.$
\end{enumerate}
\end{proposition}

\selectlanguage{english}

\end{document}